\documentclass[11pt]{article}
\usepackage{amsfonts}
\usepackage{epsfig, graphicx}
\usepackage{amssymb}
\usepackage{mathtools}
\usepackage{bbm}
\mathtoolsset{showonlyrefs}
\usepackage{amstext}
\usepackage{amsmath,enumerate}
\usepackage{theorem}
\usepackage{algorithm}
\usepackage{algpseudocode}
\usepackage[left=1in,right=1in,top=1in,bottom=1in,nohead]{geometry}
\usepackage[active]{srcltx}
\usepackage[normalem]{ulem}
\allowdisplaybreaks[2] 

\newcommand{\card}[1]{\left|#1\right|}
\newcommand{\Po}{\operatorname{Po}}
\newcommand{\Be}{\operatorname{Be}}

\newcommand\ER{Erd\H os--R\'enyi}

\newcommand{\mH}{m}
\newcommand{\vH}{v}
\newcommand{\wH}{w}
\newcommand{\WH}{W}
\newcommand{\Wk}{W}
\newcommand{\lH}{\l}

\newcommand{\dH}{d}
\newcommand{{\aH}}{a}

\def\b{\beta}

\def\D{\Delta}

\def\la{\lambda}

\def\l{\lambda}

\newcommand{\set}[1]{\left\{#1\right\}}

\newcommand{\proofstart}{{\bf Proof\hspace{2em}}}
\newcommand{\proofend}{\hspace*{\fill}\mbox{$\Box$}}

\def\E{{\sf E}}
\def\Pr{{\sf P}}

\newcommand{\ignore}[1]{}

\newtheorem{theorem}{Theorem}

\newtheorem{lemma}[theorem]{Lemma}

\newtheorem{remark}[theorem]{Remark}
\newtheorem{claim}[theorem]{Claim}
\newcommand{\den}{\mathrm{den}}
\newcommand{\aut}{\mathrm{aut}}

\newcommand{\brac}[1]{\left(#1\right)}
\newcommand{\bfrac}[2]{\brac{\frac{#1}{#2}}}

\newcommand{\beq}[2]{\begin{equation}\label{#1}{#2}\eeq}

\newcommand{\eeq}{\end{equation}}
\newcommand{\blem}[1]{\begin{lemma}\label{#1}}
\newcommand{\elem}{\end{lemma}}
\newcommand{\bth}[1]{\begin{theorem}\label{#1}}
\newcommand{\enth}{\end{theorem}}
\newcommand{\brem}[1]{\begin{remark}\label{#1}}
\newcommand{\erem}{\end{remark}}

\usepackage[usenames,dvipsnames]{xcolor}
\usepackage{xcolor}

\newcommand{\ee}{m}
\newcommand{\tvd}{\operatorname{TVD}}
\newcommand{\pa}{p}
\newcommand{\mm}{{m_\ell}}
\newcommand{\lfrac}[2]{\left. {#1} \right/ {#2}}

\newcommand{\mex}{{\binom k 2}}
\newcommand{\ai}{\alpha}
\newcommand{\bi}{\beta}
\newcommand{\FF}{F}
\renewcommand{\approx}{\sim}

\def\parens(#1){\left( {#1} \right)}

\newcommand{\ff}[2]{(#1)_{#2}}

\newcommand{\wz}{n^{-(1-\al)(1/d)-(\al)(1/d')}}
\newcommand{\al}{\alpha}
\renewcommand\ER{Erd\H os--R\'enyi}
\newcommand\GG{\mathcal{G}}
\newcommand{\muhat}{\hat{\mu}}

\newcommand{\Ebar}{{\bar E}}
\newcommand{\cond}{\mid}
\newcommand{\tabrow}[8]{#1 & #2 & #3 & #4 & #5 & #6 & #7 & #8}
\newcommand{\pab}{p_{\ai\bi}}
\newcommand{\floor}[1]{\left\lfloor #1 \right\rfloor}
\newcommand{\Fp}{F^+}
\newcommand{\Fm}{F^-}

\author{Alan Frieze\thanks{Department of Mathematical Sciences, Carnegie Mellon University. Research supported in part by NSF Grants DMS1362785, CCF1522984 and a grant(333329) from the Simons Foundation. alan@random.math.cmu.edu.},
Wesley Pegden\thanks{Department of Mathematical Sciences, Carnegie Mellon University. Research supported in part by NSF Grant DMS1363136 and the Sloane Foundation. wes@math.cmu.edu.},
Gregory Sorkin\thanks{Department of Mathematics, London School of Economics and Political Science.
g.b.sorkin@lse.ac.uk}}

\date{\today}

\begin{document}
\title{The distribution of minimum-weight cliques and
          \\ other subgraphs in graphs with random edge weights}

\definecolor{darkgray}{rgb}{0.1,0.1,0.4}
\maketitle
\begin{abstract}
We determine, asymptotically in $n$, the distribution
and mean of the weight of a
minimum-weight $k$-clique
(or any strictly balanced graph $H$)
in a complete graph $K_n$
whose edge weights are independent random values
drawn from the uniform distribution or other continuous distributions.
For the clique, we also provide explicit (non-asymptotic) bounds
on the distribution's CDF
in a form obtained directly from the Stein-Chen method,
and in a looser but simpler form.
The direct form extends to other subgraphs and other edge-weight distributions.
We illustrate the clique results for various values of $k$ and $n$.
The results may be applied to evaluate
whether an observed minimum-weight
copy of a graph $H$ in a network provides statistical evidence that
the network's edge weights are not independently distributed
but have some structure.
\end{abstract}

\section{Introduction}
In this note we consider the distribution of the minimum weight copy of a
fixed subgraph $H$ in a randomly edge weighted complete graph $K_n$;
a natural special case that may be of particular interest is when $H$ is a $k$-clique.
This seemingly natural problem seems to be absent from the literature thus far.
It can be viewed as a fixed-size version of the NP-complete Maximum Weighted Clique Problem (MWCP);
the review article \cite{Pardalos} on Maximum Clique
includes discussion of MWCP algorithms (Section 5.3)
and the performance of heuristics on random graphs (in Section 6.5).
The same article reviews applications, including in telecommunications,
fault diagnosis, and computer vision and pattern recognition.
Wikipedia's article \cite{WikiClique} on the Clique problem includes applications in
chemistry, bioinformatics, and social networks,
and many of these seem more naturally modelled as weighted than unweighted problems.
Research on fast algorithms for MWCP in ``massive graphs'' arising in practice is ongoing;
a recent example, \cite{IJCAI16}, notes applications in telecommunications and biology
(specifically, the study \cite{Moscato} of finding a 5-protein gene marker for Alzheimer's disease).

The weight distribution of a minimum-weight subgraph $H$ of a randomly weighted network
is a natural question in discrete mathematics and probability,
and in addition has statistical ramifications that
we think may be of importance in applied settings.
Taking the clique for purposes of discussion,
the distribution we seek would allow one to
judge whether the weight of the smallest-weight clique in a given network
provides statistical evidence that the network's edge weights are
not independently, uniformly random.
Because uniformity is not essential to our analysis (see Section \ref{s.X}),
an anomalous weight is evidence that the edge weights
are not independent, i.e., that the network's weights have some structure.
The statistics could alternatively be performed by repeated simulations of random networks
matching the null hypothesis
but our approach is preferable
for the usual reasons that simulation is cumbersome, it does not provide rigorous results,
and the number of simulations must be more than the reciprocal of the
desired significance level (potentially quite large),
and additionally because finding a minimum-weight clique in a graph is an NP-hard problem
that each simulation would have to solve.

We focus on the case where
each edge $e$ in $G$ has an independent weight $X_e$ which is uniform in $[0,1]$,
and $H$ is a complete graph $K_k$,
and here we derive explicit (non-asymptotic) bounds.
We also obtain asymptotic results for other graphs $H$,
and we extend both sets of results to other distributions.
We only consider finding subgraphs $H$ of the complete graph $K_n$,
but the same approach may be applicable to other networks of interest.

The \emph{density} of a graph $H$ is defined as $\den(H)=e(H)/v(H)$ where
$e(H)$ and $v(H)$ denote the number of edges and vertices of $H$,
and $H$ is \emph{strictly balanced} if
\beq{sub}{
\den(H)>\den(H')\text{ for all strict subgraphs }H' \subset H.
}

\begin{theorem}\label{thmH}
Let $H$ be a fixed strictly balanced graph with $\vH$ vertices, $\mH$ edges, $\aH$ automorphisms,
and density $\dH=\mH/\vH$.
Let the edges of $K_n$ be given independent uniform $[0,1]$ edge weights,
and let $\WH$ denote the minimum weight of a subgraph isomorphic to $H$.
Then, for any non-negative $z=z(n)$
asymptotically in $n$,
\begin{align}
\Pr\brac{\WH\geq \frac{z}{n^{1/\dH}}}
 &=\exp\set{-\frac{z^{\mH}}{\mH!\,\aH}}+o(1)
 \label{Hasymp} .
\intertext{Also,}
\E(\WH)
 &\approx \muhat := n^{-1/\dH} \frac{(\mH
   ! \, \aH)^{1/\mH}}{\mH}
     \Gamma\brac{\frac{1}{\mH}} .
 \label{0H}
\end{align}
\end{theorem}
Here $f(n) \approx g(n)$ means that $f(n)/g(n) \to 1$ as $n\to\infty$,
which we may also write as $f(n) = g(n) \, (1+o(1))$.

In the case when $H$ is a clique we have made some effort to control the $o(1)$ error, to demonstrate that our approach is useful for reasonable problem sizes.

\begin{theorem} \label{th1}
Fix $k \geq 3$, let the edges of $K_n$ be given independent uniform $[0,1]$ edge weights,
and let $\WH$ denote the minimum weight of a clique $K_k$.
Let
$$
w=\frac{z}{n^{2/(k-1)}}
  \text{ and }
\l=\binom{n}{k}\frac{w^\mex}{\mex!}\approx \frac{z^\mex}{k!\mex!} \ \
  \text{(as in \eqref{1x})}.
$$
Then, for $z \leq \min \set{ n^{2/\mex},n^{2/(k-1)} \exp(-\frac{k-1}{k-2})}$
(equivalently, $w \leq \min \set{n^{-2/k} , \exp(-\frac{k-1}{k-2}) }$),
\begin{align}
 \card{ \Pr\brac{\Wk\geq w =\frac{z}{n^{2/(k-1)}}} - \exp(-\la) }
 & \leq \, \frac87 \, \frac{(k-2)}{\binom k 2! (k-1)!^2} \:
  \: \frac{z^{\binom k 2+k-1}}{n} .
   \label{th1.2}
\end{align}
\end{theorem}

Theorem \ref{th1x} derives tighter bounds for the clique than those of Theorem \ref{th1},
but is presented later because it requires additional notation.

The proofs use the Stein-Chen method.
Section \ref{SCclique} presents the method as applied to finding a low-weight clique,
and establishes the explicit probability bounds of Theorem \ref{th1x}.
Section \ref{th1xx} outlines how explicit bounds could be obtained
for other distributions, and other subgraphs $H$.
Section \ref{bounds} simplifies (but loosens)
the bounds to give Theorem \ref{th1x}.

Section \ref{s.H} applies the Stein-Chen method to strictly balanced graphs $H$
to obtain the asymptotic probability bounds and asymptotic expectation
of Theorem \ref{thmH}.
Section \ref{s.X} generalizes Theorem \ref{thmH}
to non-uniform edge weight distributions,
as Theorem \ref{thmX},

Section \ref{sample} presents plots and tables for the application of Theorem \ref{th1x}
to various clique sizes $k$ and graph sizes $n$,
giving an indication of where our results are effective and where they need improvement.
The Conclusions in Section 8 recapitulate the results achieved
and discuss where they might be applied and how they might be extended.

\section{Stein-Chen Method} \label{SCclique}
We use a version of the Stein-Chen method given in \cite[Theorem 1]{AGG}.
We will summarize the general formulation,
which may become clearer when we then show how it applies in our case.
The formulation characterizes the distribution of a value
$X=\sum_{\ai \in I} X_\ai $
where $I$ is an arbitrary index set
and each $X_\ai$ is a Bernoulli random variable,
$X_\ai \sim \Be(\pa) , $
that is, $X_\ai$ is 1 with probability $\pa$ and 0 otherwise.
For each $z \in I$ there is a ``neighborhood of dependence''
$B_\ai \subseteq I$ with the property that
$X_\ai$ is independent of all the $X_\bi$ for $\b$ outside of $B_\ai$.
With $\l = \E X$, $Z \sim \Po(\l)$
a Poisson random variable with mean $\l$,
and
\begin{align}
 b_1 & := \sum_{\ai \in I} \sum_{\bi \in B_\ai} p_\ai p_\bi ,\label{b1b1}
 \\
 b_2 & := \sum_{\ai \in I} \sum_{\ai \neq \bi \in B_\ai} \pab
 \text{, where } \pab = \E(X_\ai X_\bi) ,
 \label{b2gen}
\end{align}
the conclusion of \cite[Theorem 1]{AGG} is that
\beq{tvd}{\tvd(X,Z) \leq  b_1+b_2 ,}
where $\tvd(X,Z)$ denotes the total variation distance between the two distributions,
i.e., the maximum, over events $E$,
of the difference in the probabilities of $E$ under the two distributions.
(The full theorem involves an additional term $b_3$ if there are weak dependencies,
but we do not need this.
Also, we have adjusted for our use of the standard definition of $\tvd$,
where \cite{AGG} defines $\tvd$ as twice this.)

In our case, with $G=K_n$, we are interested in
\beq{Xdef}{X=\text{number of $k$-cliques of $G$ weighing $\leq w$.}}

We are specifically interested in $\Pr(X=0)$ i.e. the probability that there is no such clique. We have
$$X=\sum_{\ai \in I} X_\ai $$
where the index set $I$ is the set of all $k$-cliques of $G$,
\beq{defI}{
I=\binom{[n]}{k}
}
and
\beq{defXa}{
X_\ai=
 \begin{cases}
   1, & \mbox{if clique $\ai$ has weight $\leq w$}
   \\
   0, & \mbox{otherwise} .
 \end{cases}
}
Denoting the number of edges in a $k$-clique by
\beq{defm}{\ee = \binom k 2 ,}
each  $X_\ai$ satisfies
$$ X_\ai \sim \Be(\pa) $$
where, for $w \leq 1$,
\beq{pa}{
\pa = p(w) := E(X_\ai)
 = \frac{w^\ee}{\ee!} ,
}

the probability that the sum of $m$ i.i.d.\ uniform $[0,1]$ variables is at most $w$.
The sum of i.i.d.\ uniform random variables has an Irwin-Hall distribution,
whose CDF is well known
(more on this in Section \ref{unifx}),
but it is not hard to see that when $w \leq 1$ this probability is given by \eqref{pa}.
For a clique to have total weight $\leq w$,
first, each of its $\ee$ edges must have weight $\leq w$.
With i.i.d.\ uniform $[0,1]$ random edge weights,
conditioned on the vector of $\ee$ edge weights lying in $[0,w]^\ee$,
the vector is a uniformly random point in this $\ee$-dimensional cube,
the event that the sum of its coordinates is $\leq w$
is the event that the point lies in a standard $\ee$-dimensional simplex scaled by $w$,
and the volume of this simplex is $w^\ee/\ee!$.

Immediately,
\begin{align}
\l &= \l(w) :=  \E X = \binom n k p
 = \binom{n}{k} \frac{w^\ee}{\ee!}
 \in \frac{n^k}{k!} \frac{w^\ee}{\ee!} \left[1-\frac{k^2}{2n},1\right]
 \label{1x} .
\end{align}

The neighborhood of dependence $B_\ai$ is the set of cliques that share an edge with clique $\ai$.
Breaking this down more finely, into cliques sharing $\ell$ vertices with $\ai$, we see that
\begin{align}
\card{B_\ai} &= \sum_{\ell=2}^{k} u_\ell
    \text{ where }u_\ell=\binom{k}{\ell}\binom{n-k}{k-\ell} .
    \label{Bcard}
\end{align}
Then,
\begin{align}
b_1&=\binom{n}{k} \sum_{\ell=2}^{k}u_\ell \pa^2 .
        \label{b1}
\end{align}

To calculate $b_2$, suppose that two cliques
(or indeed copies of any graph $H$)
$\ai$ and $\bi$ have $a$ edges in common and
each has $b$ edges not in the other (with $a+b=m$).
The event that $X_\al X_\b=1$, i.e., that both cliques weigh $\leq w$,
is equivalent to the shared edges weighing some $W_a \leq w$, and
both sets of unshared edges weighing $\leq w-W_a$.
As in \eqref{pa},
the CDF for the sum of
$s$ edges to be at most $r$ is $F_s(r)=r^s/s!$,
so the corresponding density is $f_s(r)=r^{s-1}/(s-1)!$.
Thus,
\begin{align}
\pab =
\pab(\al,\b,w)
 &= \int_0^w f_a(w_a) (F_b(w-w_a))^2 \, d {w_a}		 \label{pabintegral}
 \\ &= \int_0^w \frac1{(a-1)!} (w_a)^{a-1}
     \left( \frac1{b!} (w-w_a)^{b} \right)^2 \, d {w_a}   \label{integral1}
\intertext{which by a change of variable to $t=w_a/w$ is}
 &= \frac{w^{a+2b}}{(a-1)! \, (b!)^2} \int_0^1 t^{a-1} (1-t)^{2b} \, dt .
\intertext{By Euler's integral of the first kind,
 $B(x,y) = \int_{0}^1 t^{x-1} (1-t)^{y-1}dt
 = \lfrac{\Gamma(x)\Gamma(y)}{\Gamma(x+y)}$,
 this is}
 &= \frac{w^{a+2b}}{(a-1)! \, (b!)^2} \cdot \frac{(a-1)! (2b)!}{(a+2b)!}  \label{integral2}
 \\ &= \frac{w^{a+2b}}{(a+2b)!} \binom{2b}{b} .
 \label{b2exact}
\end{align}

Cliques $\ai$ and $\bi$ sharing $\ell$ vertices share $a=\binom \ell 2$ edges,
while the number of edges unique to $\bi$ is
\begin{align} \label{mm}
 b = \mm := \binom{k}{2}-\binom{\ell}{2} .
\end{align}
Thus, from \eqref{b2exact}, with \eqref{b2gen} and \eqref{Bcard},
\begin{align}
b_2
 &= \binom{n}{k} \sum_{\ell=2}^{k-1} u_\ell
 \, \pab \big( \tbinom \ell 2, \tbinom k 2-\tbinom \ell 2, w \big)   \label{b2x}
\\ &= \binom{n}{k} \sum_{\ell=2}^{k-1} u_\ell
  \, \frac{w^{m+\mm}}{(m+\mm)!} \binom{2\mm}{\mm} .
\label{b2}
\end{align}

With
\beq{Zdef}{Z \sim \Po(\l)}
the conclusion from \cite[Theorem 1]{AGG}, per \eqref{tvd}, is that
\begin{align}\label{XZ}
  \card{\Pr(W \geq w)- \exp(-\la)}
    &= \card{\Pr(X=0)-\Pr(Z=0)}
    \leq \tvd(X,Z)
    \leq b_1+b_2.
\end{align}
We have thus proved the following theorem:

\begin{theorem} \label{th1x}
Fix $k \geq 3$, let the edges of $K_n$ be given independent uniform $[0,1]$ edge weights,
and let $\WH$ denote the minimum weight of a clique $K_k$.
Let $w=z/{n^{2/(k-1)}}$,
$\l=\binom{n}{k}\frac{w^\mex}{\mex!}\approx \frac{z^\mex}{k!\mex!}$
(as in \eqref{1x}),
and $b_1$ and $b_2$ be as given in \eqref{b1} and \eqref{b2}
(calling in turn on \eqref{pa} and \eqref{Bcard}).
Then, for $z \leq n^{2/(k-1)}$ (equivalently, $w \leq 1$),
\begin{align}
 \card{ \Pr\brac{\Wk\geq w =\frac{z}{n^{2/(k-1)}}} - \exp(-\la)}
  & \leq b_1+b_2.   \label{th1.1}
\end{align}
\end{theorem}

\section{Extensions of Theorem \ref{th1x}} \label{th1xx}
In this section, we illustrate how we could, with the aid of a computer,
extend Theorem \ref{th1x} and obtain precise values for $b_1,b_2$
in more general circumstances.

Given any edge weight distribution, we can at least in principle know
the CDF $F_m(w)$ for a set of $m$ edges to have total weight at most $w$,
and the corresponding PDF $f_m(w)$.
We may then generalise \eqref{pa} to $p=F_m(w)$,
where as usual $m=\binom k 2$,
as before define $\la$ by \eqref{1x} and $b_1$ by \eqref{b1},
and compute $b_2$ through \eqref{b2x},
with $\pab$ given by \eqref{pabintegral}.

\subsection{Uniform edge weights} \label{unifx}
With edge weights uniformly distributed as before, but removing Theorem \ref{th1x}'s
restriction to $w \leq 1$,
the sum of $m$ uniform weights has Irwin-Hall distribution,
with distribution and density functions
\begin{align}\label{Fnx}
  F_m(w) &= \frac1{m!} \sum_{i=0}^{\floor w} (-1)^i \binom m i (w-i)^m
  \\
  f_m(w) &= \frac1{(m-1)!} \sum_{i=0}^{\floor w} (-1)^i \binom m i (w-i)^{m-1} ;
\end{align}
see \cite[eq(12)]{Hall}, \cite{Irwin}, \cite[eq(26.48)]{Johnson2}.

Given $k$ (thus $m$) and $w$ it is straightforward to calculate $F_m(w)$ and
thus $p$, $\la$, and $b_1$.
Calculating $b_2$ reduces, for each $\ell$ in \eqref{b2x},
to calculating $\pab$ through the integral in \eqref{pabintegral}.
Break the range of integration into intervals within which neither $w$ nor $w_a$
takes an integral value,
by splitting at the points where either does take an integral value;
since $w \leq m$ there are at most $2m$ such points.
Within each subintegral,
the integrand is the product of a polynomial in $w_a$ (from $f_a(w_a)$)
and a polynomial in $(w-w_a)$ (from the square of $F_b(w-w_a)$).
Each term of this can be integrated as an Euler integral of the first kind,
as was done in going from \eqref{integral1} to \eqref{integral2},
or indeed expanded to a polynomial in $w_a$
(the powers are all bounded in terms of $k$)
where each term can be integrated straightforwardly
(even as an indefinite integral).

In principle, then, we can extend Theorem \eqref{th1x} to all $w \leq m$.
Indeed, for each $w$ we can produce $p$ and $\pab$,
thus giving $b_1$ and $b_2$ as explicit functions of $n$.
(We cannot get an explicit function of $w$, at least by the method above,
because the partition of the integral into sub-integrals is
different for each $w$.)
In practice, a naive implementation in Maple struggles with $k=10$
both in computation time and numerical stability.

\subsection{Exponential edge weights} \label{expx}
We may also consider exponentially distributed edge weights.
Without loss of generality we assume rate 1; anything else is a simple rescaling.
The sum of $m$ rate-1 exponentials has Erlang distribution with
well known distribution and density functions
that can be stated in a variety of forms including
\begin{align}\label{Erlang}
  F_m(w) &= 1-  e^{-w} \, \sum_{i=0}^{m-1} \frac{w^i}{i!}
  \\
  f_m(w) &= \frac{w^{m-1} e^{-w}}{(m-1)!}
  ;
\end{align}
see for example \cite{WikiErlang} and \cite[eq(17.2)]{Johnson1}.

Again given $w$ it is straightforward to compute $p$, $\la$, and $b_1$,
while computing $b_2$ requires computing $\pab$ for each $\ell$
in the sum in \eqref{b2x} and the key is to evaluate \eqref{pabintegral}.
In this case the form of $F_b(w)$ means that the integrand is a finite sum
(with length a function of $k$),
each term of which has form a constant (with respect to $w$) times
$(w_a)^a e^{-w_a}$ (coming from $f_a(w_a)$) times
$(w-w_a)^r e^{-s (w-w_a)}$ (from the square of $F_b(w-w_a)$),
for some integers $r$ and $s$.
These may in turn be expanded to terms of form
$(w_a)^r e^{s \, w_a}$.
Each of these is integrable (even as an indefinite integral);
alternatively, each definite integral (over $w_a$ from $0$ to $w$)
is a lower incomplete gamma function.

It thus seems feasible, if unenviable, to extend Theorem \ref{th1x}
to provide explicit bounds for exponential edge weights.
In practice, Maple has little difficulty with the calculations through $k=6$
but they quickly get more difficult:
at $k=10$ is is challenging to calculate
even a single $\ell$ the corresponding function $\pab(w)$.

\subsection{Subgraphs $H$ other than cliques} \label{exactH}
Generalising Theorem \ref{th1x} to subgraphs $H$ other than cliques
appears straightforward.
The neighborhood of dependence of a given copy of $H$ needs more careful treatment,
but this can be done in this non-asymptotic setting precisely as presented
in Section \ref{s.H} for asymptotic calculations.
The explicit calculations here do not even require that $H$ be strictly balanced,
but it can be expected that the error term $b_2$ will be large if it is not
(for the same reasons that the strict balance is generally required in application
of the second moment method).

\section{Calculating bounds} \label{bounds}
Given values of $k$ and $w$, in practice one would apply Theorem \ref{th1x}
using \eqref{th1.1}, calculating $b_1+b_2$ exactly from \eqref{b1} and \eqref{b2},
as indeed we do in Section \ref{sample}.
However, to characterize the quality of the estimate of $\Pr(W \geq w)$
we derive an upper bound on $b_1+b_2$ as a relatively simple (summation-free) function of
$k$ and $w$.

We start with $b_2$, the more difficult and (as we will see)
larger of these two parameters.
First, in lieu of \eqref{b2exact},
we observe that for clique $\bi$ to have weight at most $w$,
the $\mm$ edges unique to it must have total weight $\leq w$,
and therefore
\begin{align}
\Pr(X_\bi =1 \mid X_\ai =1) \leq \frac{w^\mm}{\mm!}
\leq \frac{w^{\mm}}{(k-1)!}.                        \label{prweak}
\end{align}
In the first inequality we have used that since these edges are unique to $\b$,
conditioning on the weight of $\ai$ being at most $w$ is irrelevant.
The second  is simply because, over the range of $\ell$ from 2 to $k-1$,
$\mm \geq \binom k 2 - \binom{k-1}2 = k-1$.
Also,
\begin{align}\label{ubound}
  u_\ell & \leq \binom k \ell \, \frac{n^{k-\ell}}{(k-\ell)!} .
\end{align}
Substituting \eqref{prweak} and \eqref{ubound} into \eqref{b2x}, it follows that
\begin{align} \label{b2upper}
b_2 & \leq b_2' := \binom n k  \sum_{\ell=2}^{k-1} v_\ell
 \text{ where }
 v_\ell = \binom k \ell \, \frac{n^{k-\ell}}{(k-\ell)!}\, \pa \, \frac{w^\mm}{(k-1)!}   .
\end{align}

\begin{claim} \label{vconvexclaim}
Assuming that $w \leq \min \set{ n^{-2/k}, \exp(-\frac{k-1}{k-2}) }$,
over $2 \leq l \leq k-1$, $v_\ell$ is maximized at  $\ell=k-1$.
\end{claim}

\proofstart
We first claim that
$v_\ell$ is log-convex over $2 \leq \ell \leq k-1$,
so that the maximum occurs either at $\ell=2$ or $\ell=k-1$.
For $k=3$ and 4 this is trivial. Otherwise,
for $3\leq \ell\leq k-2$,
\beq{eq3a}{
\frac{v_{\ell+1}}{v_\ell}=\frac{(k-\ell)^2}
{w^\ell n (\ell+1)}
}
and
\begin{align}
\frac{v_\ell^2}{v_{\ell-1}v_{\ell+1}}
 &=\frac{\ell+1}{\ell}\cdot \bfrac{k-\ell+1}{k-\ell}^2  w.
\label{vconvex1}
\end{align}
To establish that $v_\ell$ is log-convex
over $2 \leq l \leq k-1$
it suffices to show that
\eqref{vconvex1} is $\leq 1$
over $3 \leq l \leq k-2$.
Using $1+x \leq \exp(x)$,
from \eqref{vconvex1} we have
\begin{align}
\frac{v_\ell^2}{v_{\ell-1}v_{\ell+1}}
&\leq w\exp\set{\frac 1 \ell + \frac 2 {k-\ell}}
 \leq w \exp\set{\frac{k-1}{k-2}}
 \leq 1
  \label{xx1} ,
\end{align}
where the final inequality is by hypothesis
and the previous one because $\frac 1 \ell + \frac 2 {k-\ell}$
is convex, so its maximum occurs at one of the extremes,
either $\ell=3$ or (in fact) $\ell=k-2$.

Thus, $v_\ell$ is log-convex and its maximum occurs either at
$v_2$ or $v_{k-1}$.
However,
$$
\frac{v_2}{v_{k-1}} =
 \frac{\binom k 2 n^{k-2} w^{\binom k 2-\binom 2 2} / (k-2)!}
 {\binom k {k-1} n^1 w^{\binom k 2-\binom {k-1} 2} / (1!)}
= \frac{k-1}{2(k-2)!} (n w^{\frac12 k})^{k-3}
\leq 1,
$$
by the hypothesis that $w \leq n^{-2/k}$.
Thus, $v_2 \leq v_{k-1}$, proving the claim.
\proofend

It follows from \eqref{b2upper} and Claim \ref{vconvexclaim}  that
\begin{align}
b_2 \leq b_2'
 & \leq \binom n k (k-2) \cdot v_{k-1}
 \leq
 \binom n k (k-2) \cdot
 \binom k {k-1} \, \frac{n^1}{1!} \, \frac{w^{\binom k 2}}{\binom k 2!} \, \frac{w^{k-1}}{(k-1)!}\\
& \leq
 \frac{(k-2)}{\binom k 2! (k-1)!^2} \: n^{k+1} w^{\binom k 2+k-1} .
\label{eq3}
\end{align}
Recalling \eqref{b1}, using \eqref{ubound}, and by analogy with \eqref{b2upper},
\begin{align} \label{b1upper}
b_1 & \leq b_1' := \binom n k  \sum_{\ell=2}^{k} v'_\ell
 \text{ where }
 v'_\ell = \binom k \ell \, \frac{n^{k-\ell}}{(k-\ell)!} \, \pa^2 .
\end{align}
From the definitions of $v_\ell$ and $v'_\ell$ in \eqref{b2upper} and \eqref{b1upper},
for $2 \leq \ell \leq k-1$,
$$
\frac{v'_\ell}{v_\ell }
 = \frac{p} {\lfrac{w^\mm}{(k-1)!}}
  = \frac {\lfrac{w^m}{m!}} {\lfrac{w^\mm}{(k-1)!}}
  \leq \frac13 w^{\binom \ell 2}
  \leq \frac13 w\leq \frac13\exp\brac{-\frac{k-1}{k-2}}\leq \frac{1}{3e} .
$$
It follows (referring now to the definitions of $b'_2$ and $b'_1$ in \eqref{b2upper} and \eqref{b1upper})
that the sum of all but the $\ell=k$ term in $b'_1$, which has no counterpart in $b'_2$,
is at most $\frac1{3e} b'_2$.
Also, from \eqref{b1upper}, the $\ell=k$ term of $b'_1$ is
is small compared with its $\ell=k-1$ term:
$$
\frac{v'_k}{v'_{k-1}}
 = \frac{\binom k k n^{0} / 0!}{\binom k{k-1} n^1/1! }
 = \frac1{k n}
 \leq \frac1{9} .
$$
That is, the last term in the summation for $b'_1$ is at most $1/9$ times the second-last,
therefore it is at most $1/9$ times the sum of all but the last,
so that
\begin{align} \label{b1bound}
 b_1 \leq b'_1 \leq \brac{1+\frac1{9}} \frac1{3e} b'_2
 < \frac{b'_2}{7} .
\end{align}
Equation \eqref{th1.2} of Theorem \ref{th1}
follows from
\eqref{b1bound}, \eqref{eq3}, and \eqref{th1.1} of Theorem \ref{th1x},
since substituting $w = z/ n^{2/(k-1)}$ into the term
$n^{k+1} w^{\binom k 2+k-1}$ of \eqref{eq3}
gives $z^{\binom k 2 +k-1}/{n}$.

\section{Strictly balanced $H$} \label{s.H}
\newcommand{\nksets}{\binom{[n]}v}

Let $H$ be a strictly balanced graph with $\vH$ vertices, $\mH$ edges,
and automorphism group $\aut(H)$ of cardinality $\aH$.
We apply the Stein-Chen method in parallel with Section \ref{SCclique}.
A bit more care is needed in defining the index set $I$.
Think of a copy of $H$ in $G$ as defined by a set
of $\vH$ vertices of $G$,
together with a 1-to-1 mapping
from these vertices to those of $H$.
The set $S$ of vertices is drawn from the collection $\nksets$,
the set of all $\vH$-element subsets of $[n]$, with $\card{\nksets} = \binom n \vH$.
Taking the elements of $S$ in lexicographic order, the mapping into $V(H)$ is given by
a permutation $\pi$ of the values $1 \ldots \vH$, taken modulo the automorphism group of $H$.
Thus we may draw $\pi$
from a set $L$ of permutations, with
$$ \card L = \frac{\vH!}{\aH}  , $$
$L$ consisting of one permutation from each equivalence class.

Then, in analogy with Section \ref{SCclique}'s equation \eqref{defI}, here
\beq{defIH}{
I=\nksets \times L .
}
With no change from before, the number of copies of $H$ of weight $\leq w$
is given by a random variable $X = \sum_{\ai \in I} X_\ai$,
the $X_\ai$ Bernoulli random variables.
In parallel with \eqref{pa}, each $X_\ai$ has expectation
\beq{EXH}{
\pa := E(X_\ai) = \frac{w^\ee}{\ee!} ,
}
assuming $w \leq 1$, and in parallel with \eqref{1x},
\beq{1xH}{
\l = \l(w) := \E X = \card I \cdot \pa = \binom n \vH \frac{\vH!}{\aH} \frac{w^\ee}{\ee!} .
}

As before, we focus on $b_2$ and then treat $b_1$.
The structure of dependent events here is a bit subtle, and an example may be useful.
Suppose $H$ is the 2-path with edges $\set{1,2}$ and $\set{2,3}$.
Its only automorphism is relabeling 123 to 321.
On vertices, say $S=\set{9, 11, 15}$ of $G$, there are then 3 index sets:
the set $\set{9,11,15}$ in combination with any permutation chosen from
$L=\set{123, 231, 312}$,
corresponding respectively to paths 9--11--15, 11--15--9, and 15--9--11.
(For instance in 231, we take the 2nd, 3rd, and 1st elements of the set in that order.)
The permutations 321, 132, and 213 are eliminated from $L$ by automorphism,
and correspond respectively to paths 15--11--9, 9--15--11, and 11--9--15 already listed.
Suppose $\ai$ is vertex set $\set{9,11,15}$ with permutation 123,
giving path 9--11--15.
Consider all possible dependent indices $\bi$ on vertices $\set{11,15,18}$.
With $\pi=$ 123,
$\bi$ gives path 11--15--18, sharing edge 11--15 with $\ai$.
With $\pi=$ 312,
$\bi$ gives path 18--11--15, again sharing edge 11--15 with $\ai$.
Finally, with $\pi=$ 231,
$\bi$ gives path 15--18--11, sharing no edges with $\ai$,
and thus $\bi \notin B(\ai)$.

For $\bi \in B(\ai)$, the pair $(\ai,\bi)$ describes an overlapping pair of copies of $H$:
a pair of labeled graphs $(V_1,E_1)$ and $(V_2,E_2)$, each isomorphic to $H$,
$V_1,V_2 \in \nksets$,
with $\card{E_1 \cap E_2} \geq 1$ and
thus $\ell:=\card{V_1 \cap V_2} \geq 2$.
Their union is a graph $\FF = (V_1 \cup V_2, E_1 \cup E_2)$,
and for strictly balanced $H$, $\den(\FF) > \den(H)$.
(This is easy to show, well known, and the reason for introducing strict balance;
one early reference is \cite[eq(4.3)]{Bol81}.
Any such graph $\FF$ has at most $2\vH-2$ vertices,
so up to the labeling of the vertices there are only finitely many possibilities.
Let
\begin{align} \label{d'def}
d'=d'(H)
\end{align}
be the minimum density of all such graphs $\FF$.
For example, if $H=K_k$ then $d'$ is obtained for two $k$-cliques sharing $k-1$ vertices,
and
$d'=\frac{2\binom k 2-\binom{k-1}2}{2k-(k-1)}=\frac{(k-1)(k+2)}{2(k+1)}>d=\frac{k-1}{2}$.

For copies $\ai$ and $\bi$ both to have weight $\leq w$,
their union graph $\FF$ must have weight $\leq 2w$,
and as in \eqref{pa} this event has probability $(2w)^{m(\FF)}/m(\FF)!$,
assuming $2w \leq 1$.
Thus, if the two copies overlap in $\ell$ vertices,
implying that $m(\FF) \geq d'(2\vH-\ell)$, we have
\begin{align}
 \E(X_\ai X_\bi)
  & \leq \frac{(2w)^{m(\FF)}}{m(\FF)!}
   \leq \frac {(2w)^{d'(2\vH-\ell)}} {(d'(2\vH-\ell))!} .
 \label{PabH}
\end{align}

It follows that
\begin{align}
 b_2
  &\leq \card I \sum_{\ell=2}^\vH \binom \vH \ell  \binom{n-\vH}{\vH-\ell} \card L
        \frac {(2w)^{d'(2\vH-\ell)}} {(d'(2\vH-\ell))!}
  \\&= O\brac{ n^\vH \sum_{\ell=2}^\vH n^{\vH-\ell} w^{d'(2\vH-\ell)}}
  \\&= O\brac{ \sum_{\ell=2}^\vH  (n \, w^{d'})^{2\vH-\ell} }
  \\&= O\brac{ (n \, w^{d'})^\vH } = o(1),
  \label{bigOH}
\end{align}
the last pair of inequalities holding subject to the condition that
\begin{align} \label{wcond}
w=o(n^{-1/d'}) \quad \text{or equivalently} \quad
  z = {w} \, {n^{1/d}} = o(n^{1/d-1/d'}) .
\end{align}

Note that unlike in \eqref{b2}
the sum here includes $\ell=\vH$
but nonetheless capitalizes on $\bi \neq \ai$ from the definition of $b_2$ (see \eqref{b2gen}):
the vertex sets of $\bi$ and $\ai$ may be equal but the index sets themselves are different,
so that the union graph $\FF$ is not isomorphic to $H$, and therefore $\den(\FF) \geq d' >d$.

Now compare $b_1$ and $b_2$ from their definitions in \eqref{b1b1} and \eqref{b2gen}.
For each term $(\ai,\bi)$ common to both sums,
the summand in $b_1$ is $p_\ai p_\bi = p^2 = O(w^{2m})$,
and (with $w \leq 1$) this is of smaller order than the corresponding summand in $b_2$
(see \eqref{b2gen} and \eqref{PabH}),
which is of order $w^{m(\FF)}$.
($\FF$ is formed of two copies of $H$ sharing at least one edge, thus $m(\FF) \leq 2m-1$.)
Only the terms $(\ai,\ai)$ are unique to $b_1$,
within $b_1$ they are fewer than the other terms, and all terms are equal,
so they do not change the order of $b_1$.
It follows that
$$b_1=O(b_2) . $$
We have established that, subject to \eqref{wcond},
\begin{align} \label{thm1.1nearly}
|\Pr(\WH \geq \wH) - e^{-\lH}| \leq b_1 + b_2 = O\brac{ (n \, w^{d'})^\vH } = o(1) .
\end{align}

With $z=w \, n^{1/d}$, and using the usual falling-factorial notation,
observe from \eqref{1xH} that
\beq{lap}{
 \l = \frac{\ff n v}{a} \frac{w^m}{m!}
  \approx \frac{n^\vH}{\aH} \frac{w^\ee}{\ee!} = \frac{z^m}{m! \, \aH} =: \l' .
}
By the intermediate value theorem, there is a point $\l'' \in [\l,\l']$
at which $\frac d {d\l} e^{-\l} = (\exp(-\l)-\exp(-\l'))/(\l-\l')$.
It follows that
\begin{align}\label{lambdas}
  {e^{-\l'}-e^{-\l}}
   & = { (\l'-\l) \cdot \frac{d}{d \l} \exp(-\la) \big|_{\l=\l''}}
   \\&=  o(1) \l'' \cdot \exp(-\l'')
  =o(1) ,
\end{align}
using that both $\l$ and $\l$ are $\l''(1+o(1))$ and that $\l'' \exp(-\l'') \leq 1/e$ for any $\l'' \geq 0$.
Now, in \eqref{thm1.1nearly} substitute $w=z/n^{1/d}$,
yielding
\begin{align} \label{PrHmain}
\card{\Pr\brac{\WH \geq \frac{z}{n^{1/\dH}}}-e^{-\l'}}
 &\leq
\card{\Pr\brac{\WH \geq \frac{z}{n^{1/\dH}}}-e^{-\l}} + \card{e^{-\l}-e^{-\l'}} = o(1) .
\end{align}
This completes the proof of \eqref{Hasymp}
subject to \eqref{wcond}, i.e.,
for $z=o(n^{1/d-1/d'})$.

To extend this to all $z=z(n)$, we will observe that there is a weight threshold $w_0$
(see \eqref{w0} below)
where $w_0$ is large compared with $n^{-1/d}$ so that a cheap copy of $H$
(cheaper than $w_0$) is  almost certainly present,
but $w_0$ is small compared with $n^{-1/d'}$ so that an overlapping pair of cheap copies is almost certainly
not present and thus the error bound $b_1+b_2$ is small.
Values $w<w_0$ are controlled by the previous case,
while for values $w>w_0$, \eqref{Hasymp} holds trivially because all its terms are $o(1)$.
The same thresholding around $w_0$ will be used shortly in proving \eqref{0H}.
We now implement this idea.

For any $0 < \al < 1$
(throughout, $\al=1/2$ will do),
define
\begin{align} \label{w0}
  w_0 &= \wz \quad \text{and} \quad
  z_0 = {w_0} \, {n^{1/d}} = n^{\al(1/d-1/d' )} .
\end{align}
By construction,
\begin{align} \label{w0small}
w_0 \, n^{1/d'} = n^{(1-\al)(1/d'-1/d)} = o(1)
\end{align}
so that
\eqref{wcond} is satisfied,
while at the same time
\begin{align} \label{z0big}
 z_0 = w_0 \, n^{1/d} = n^{(\al)(1/d-1/d')}
  = \omega(1) .
\end{align}

A putative counterexample to Theorem~\ref{thmH} equation \eqref{Hasymp}
consists of an infinite sequence $z=z(n)$ for which the error terms are not $o(1)$.
Divide such a sequence into two subsequences according to whether $z \leq z_0$ or $z>z_0$.
We have just established that the subsequence with $z \leq z_0$ must give error terms $o(1)$,
so consider the subsequence with $z>z_0$.
Here,
\begin{multline*}
\left|\Pr\brac{\WH\geq \frac{z}{n^{1/\dH}}}-\exp\set{-\frac{z^{\mH}}{\mH!\,\aH}}\right|
 \leq
\Pr\brac{\WH\geq \frac{z_0}{n^{1/\dH}}}+\exp\set{-\frac{z^{\mH}}{\mH!\,\aH}}
 \\ \leq
 \left( \exp\set{-\frac{z_0^m}{\mH!\,\aH}}+o(1) \right) +\exp\set{-\frac{z^{\mH}}{\mH!\,\aH}}
  =o(1),
\end{multline*}
where the application of \eqref{PrHmain} is justified by \eqref{w0small}
and both exponential terms are small because of
\eqref{z0big} and $z \geq z_0$.
This completes the proof of \eqref{Hasymp} for all $z$.

We now turn to Theorem \ref{thmH}, Equation \eqref{0H}.
We have
\begin{align}
\E(\WH)
&= \int_{0}^{\mH}\Pr(\WH\geq w) \, dw
= \int_{0}^{w_0}\Pr(\WH\geq w) \, dw
 + \int_{w_0}^{\mH}\Pr(\WH\geq w) \, dw .
\end{align}
Setting
\begin{align}
c &= \frac{\ff n v}{a \, m!}
\end{align}
and substituting \eqref{thm1.1nearly} into the first integral
(as justified by \eqref{w0small} and \eqref{wcond})
gives
\begin{align}
\E(\WH)
&=\int_{0}^{w_0}
 \parens( {e^{-c w^m} +
        O\parens({ {(n\,w^{d'})}^v }) }) \, dw
 + \int_{w_0}^{\mH} \Pr(\WH\geq w) \, dw
\\
&=
 \int_{0}^{w_0} e^{-c w^m} \, dw
 + \int_{0}^{w_0}
        O\parens({ {(n\,w^{d'})}^v }) \, dw
 + \int_{w_0}^{\mH} \Pr(\WH\geq w) \, dw
 \label{EW1}
\\
&= (1+o(1))
 n^{-1/d} \frac{(a \, m!)^{-1/m}}{m}
   \; \Gamma\parens(\frac1m)
  + o(n^{-1/d}) + O(\exp(-n^{\Omega(1)})) ;
 \label{EW2}
\end{align}
we will prove \eqref{EW2} by considering
each of the three integrals in \eqref{EW1} in turn.%
\footnote{%
Landau notation does not normally presume the sign of the quantity
in question,
but in error expressions like \eqref{term2small} and \eqref{EW2}
we mean for $\Omega$ to denote a positive quantity.
}
From \eqref{EW2},
Theorem \ref{thmH}, Equation \eqref{0H} follows immediately.

The first integral in \eqref{EW1} is the principal one.
Let $x=c w^m$ so that
$w=(x/c)^{1/m}$ and $dw=c^{-1/m} \frac1m x^{\frac1m-1} dx$.
Then
\begin{align}
\int_0^{w_0} e^{-cw^m} dw
 &= \frac{c^{-1/m}}{m} \:
  \int_0^{c w_0^m} e^{-x} \:  x^{\frac1m-1} dx
 \\& \sim n^{-1/d} \frac{(a \, m!)^{-1/m}}{m}
   \; \Gamma\parens(\frac1m) .
\end{align}
The asymptotic equality above follows from considering the two multiplicands separately.
For the first multiplicand,
$c^{-1/m}/m$,
we just observe that $c$'s term $(\ff n v)^{-1/m} \sim n^{-v/m} = n^{-1/d}$.
For the second multiplicand, the integral,
the upper limit of integration is tending to infinity:
$c \, w_0^m$ is of order
$n^v w_0^m
 = (n^{1/d} w_0)^m
 \to \infty
$,
by \eqref{z0big}.
Thus the integral is asymptotic to
$ \int_0^{\infty} e^{-x} \:  x^{\frac1m-1} dx$,
which is equal to
 $ \Gamma\parens(\frac1m) $:
it is an example of Euler's integral of the second kind,
$\Gamma(t)=\int_{0}^\infty x^{t-1}e^{-x}dx$.

For the second integral in \eqref{EW1},
\begin{align}
\int_{0}^{w_0}   O\parens({ {(n\,w_0^{d'})}^v })  \, dw
 &= O\parens({ w_0 \cdot (n^{1/d'} w_0)^{d' v} })
 \intertext{which, from $n^{1/d'} w_0 = o(1)$ by \eqref{w0small} and $v \geq 2$ is}
 &= O\parens({ w_0 \cdot (n^{1/d'} w_0)^{2 d'} })
 = O\parens({ w_0 \cdot (n \, w_0^{d'})^2 }) .
\end{align}
To show that this is $o(n^{-1/d})$ as claimed in \eqref{EW2}
means showing that, when multiplied by $n^{1/d}$, it is $o(1)$.
This follows from
\begin{align} \label{term2small}
 n^{1/d} \cdot w_0 \cdot (n \, w_0^{d'})^2
  &=
n^{(1/d-1/d') \, (2 \al d'+\al-2 d')}
  = n^{-\Omega(1)} = o(1) ,
\end{align}
the final two inequalities holding if $2 \al d'+\al-2 d'<0$,
i.e., if $\al < 2d'/(2d'+1)$.
Recall from \eqref{d'def} that $d'$ is the density of a graph $F$
describing an overlapping pair of copies of $H$;
say $F$ has $m'$ edges and $v'$ vertices.
Since $F$ is connected, $v' \leq m'+1$ and $d'=m'/v' \geq m'/(m'+1)$.
Since there is at least one edge shared between the two copies
and one edge unique to each copy, $m' \geq 3$,
so $d' \geq 3/4$.
Thus \eqref{term2small} holds
for any $\al < (2 \cdot 3/4)/(2 \cdot 3/4+1) = 3/5$.
Fixing for example $\al=1/2$
(in all other parts of the proof, any $\al$ strictly between 0 and 1 will do),
the integrated error term is indeed $o(n^{-1/d})$.

For the third integral in \eqref{EW1},
while $\Pr(W \geq w)$ of course decreases with $w$,
it is difficult for us to capitalize on this since
our estimates cannot be applied for $w>n^{1/d'}$
where condition \eqref{wcond} is violated.
If as for the second integral
we reason through $\Pr(W \geq w) \leq \Pr(W \geq w_0)$,
the estimate is not good enough:
we get an expression like \eqref{term2small}
but with its integration range of $w_0$ replaced by $\Theta(1)$,
giving $n^{-\frac1d [(d'-d)(2-2\al)-1]}$,
and if $d'$ and $d$ are nearly equal the exponent is not negative
for any $\al$ between 0 and 1.

\begin{claim} \label{PrExp}
For any $0 < \al < 1$, with $w_0$ given by \eqref{w0},
$\Pr(W > w_0) \leq \exp(-n^{\Omega(1)})$.
\end{claim}

\proofstart
First, we claim that an \ER\ random graph $G \sim \GG(n, w_0)$
contains a copy of $H$ w.p.\ $>1/2$.
This can be obtained as a classical application of the second-moment method,
but it also follows trivially from \eqref{PrHmain}:
The set of edges of weight $\leq w_0$ forms a random graph $G$,
the claim is that this subgraph includes a copy of $H$ w.p.\ $>1/2$,
and \eqref{PrHmain} says that with even higher probability (namely $1-o(1)$)
there exists such a copy with additional properties
(not only is each edge weight $\leq w_0$, but the total is also $\leq w_0$).

In particular, the set of edges of weight $\leq w_0$ forms a
random graph $G_1 \sim \GG(n, w_0)$,
and $G_1$ contains a copy of $H$ w.p. $>1/2$.

Now form a second, independent random graph $G_2 \sim \GG(n, w_0)$,
each of whose edges has weight $\leq 2 w_0$,
by the following standard trick.
For an edge appearing in $G_1$, accept it into $G_2$ w.p.\ $w_0$.
For an edge not appearing in $G_1$,
accept it into $G_2$
with probability $1-w_0$ \emph{if} its weight is between $w_0$ and $2w_0$,
and otherwise reject it.
Note that in the second case, the weight is in the range $(w_0,2w_0)$
with probability $w_0/(1-w_0)$,
and then we take it only w.p.\ $1-w_0$,
for a net probability of $w_0$.
Thus, each edge appears in $G_2$ with probability exactly $w_0$
independent of $G_1$ and the other edges.

As an \ER\ random graph, $G_2$ contains a copy of $H$ w.p.\ $>1/2$,
independently of $G_1$, and in such a copy every edge has weight $\leq 2w_0$.

Repeat this process for graphs $G_3$, \ldots, $G_k$.
With probability $\geq 1-2^{-k}$ at least one of these graphs contains
a copy of $H$,
and if so all its edges have weight $\leq k w_0$
for total weight $\leq k \, m \, w_0$.

To get the claim,
given $\al$, choose smaller constants $0 < \al'' < \al' < \al$.
These give rise to corresponding values $w_0'' < w_0 ' < w_0$,
and the ratios $\D'' = w_0'/w_0''$ and $\D'=w_0/w_0'$
are both of order $n^{\Omega(1)}$.
Given $G_1,G_2,\ldots G_k$ with $k=\D''$, we look for $H$
in the $k$ copies of $G(n,w_0'')$:
we find such a copy of $H$ w.p.\ $1-2^{-\D''} = 1-\exp(-n^{\Omega(1)})$,
and any such copy has weight at most
$k \cdot m \cdot w_0'' \leq \D'' \cdot \D' \cdot w_0'' = w_0$.
\proofend

From Claim \ref{PrExp} it is immediate that
\begin{align}
 0 <
 \int_{w_0}^{\mH} \Pr(\WH\geq w) \, dw
 &\leq
 m \cdot \Pr(\WH \geq w_0)
 = m \cdot \exp(-n^{\Omega(1)}) ,
\end{align}
and we absorb the constant $m$ into the $\Omega$.
This concludes analysis of the third integral in \eqref{EW1},
and thus concludes the proof of
Theorem \ref{thmH}, Equation \eqref{0H}.

\section{Extension of Theorem \ref{thmH}}\label{s.X}

Theorem \ref{thmH} extends to distributions other than uniform on $[0,1]$;
such extensions are common in situations where, intuitively,
only edges with very small weights are relevant.

\begin{theorem} \label{thmX}
The conclusions of Theorem \ref{thmH} hold under the same hypotheses
except that now the edge weights are i.i.d.\ copies of any
non-negative random variable $X$ with finite expectation and
a continuous distribution function $F$
that is differentiable from the right at 0,
with slope $F'(0)=1$.
\end{theorem}

The assumption that $F'(0)=1$ is without loss of generality.
As is standard, it can be extended to a variable $X$ for which $F'(0)=c$, for any $c>0$,
simply by rescaling: applying the theorem to $c \, X$.

To prove the theorem,
couple $X$ with a random variable $U=F(X)$.
As the quantile of $X$, $U$ is distributed uniformly on $[0,1]$.

As in \eqref{w0}, fix $0 < \al < 1/2$ and define $w_0$ accordingly;
recall that $w_0 \to 0$ as $n \to \infty$.
Any edge of weight $w_U \leq  2 w_0 \to 0$ in model $U$
has weight $w_X = F^{-1}(w_U)=w_U \, (1+o(1))$ in model $X$,
since $F'(w_0) \to 1$.
Symmetrically,
any edge with weight $w_X \leq 2 w_0$ in model $X$
has weight $w_U = F(w_X)=w_X \, (1+o(1))$ in model $U$.

Let $H_U$ and $H_X$
denote the lowest-weight copies of $H$ in the two models,
and $W_U$ and $W_X$ the corresponding optimal weights.
If $W_U \leq w_0$ then $W_X \leq W_U (1+o(1))$,
since $H_U$ would give such an $X$-weight
(each of its constituent edges has weight $w_U \leq w_0$,
thus asymptotically equal weight $w_X$)
and the weight of $H_X$ may be even smaller.
Taking the same hypothesis not the symmetric one as might be expected,
if $W_U \leq w_0$,
then $W_X \leq 2 \, w_0$
(this is why we introduced the factor of 2),
in which case
(now symmetrically)
$W_U \leq W_X (1+o(1))$.

Thus, if $W_U \leq w_0$ --- call this ``event $E$'' --- then
\begin{align} \label{wRatios}
W_X/L \leq W_U \leq W_X L
\end{align}
for some $L = L(w_0) = 1+o(1))$
in the limit $n \to \infty$ and thus $w_0 \to 0$.
Let $\Ebar$ be the complementary event
and recall from Claim \ref{PrExp} that $\Pr(E) = 1-\exp(n^{-\Omega(1)})$.

We now prove the distributional result \eqref{Hasymp} for $X \sim F$.
Rewrite \eqref{Hasymp} as
\begin{align} \label{Hasympf}
\Pr(W_U \geq w) &= f(w) + o(1) ,
\end{align}
where
$f(w) =
 \exp\set{-\frac{(w n^{1/d})^{\mH}}{\mH!\,\aH}}$.
Note that if we change $w$ by a factor $L=1+o(1)$
then the argument of the exponential changes by a factor $L^m=1+o(1)$
and thus, by \eqref{lambdas},
$f$ changes by an additive $o(1)$, i.e.,
$f(wL)=f(w)+o(1)$.

Given any $w$,
and conditioning on event $E$,
\begin{align}
\Pr(W_X \geq w \cond E)
 & \geq \Pr(W_U/L \geq w) \quad \text{from \eqref{wRatios}}
 \\& = \Pr(W_U \geq w \, L)
 \\ &= f(w \, L) + o(1) \quad \text{from \eqref{Hasymp} and \eqref{Hasympf}}
 \\ &= f(w) + o(1) \quad \text{by the argument below \eqref{Hasympf}} .
\end{align}
Symmetrically,
$
\Pr(W_X \geq w \cond E)
  \leq \Pr(W_U \geq w/L)
 = f(w) + o(1)
$
and thus
$
\Pr(W_X \geq w \cond E)
 = f(w) + o(1) .
$
For any event $A$, and any event $E$ of probability $1-o(1)$
it holds that
$\Pr(A) = \Pr(A \cond E) + o(1)$,
so here it follows that
$
\Pr(W_X \geq w)
 = \Pr(W_X \geq w \cond E) + o(1)
 = f(w) + o(1) .
$

This completes the proof of the distributional result.
We now prove the expectation result \eqref{0H} for $X \sim F$.

Let $w_0$ and event $E$ be as above.
Recall from Claim \ref{PrExp} that $\Pr(\Ebar) = \exp(-n^{\Omega(1)})$,
and from Theorem \ref{thmH} that $\E(W_U) = \Theta(n^{-1/d})$.
By the law of total expectation,
\begin{align}
\E(W_U)
 &= \E(W_U \cond E) \Pr(E) + \E(W_U \cond \Ebar) \Pr(\Ebar)
\end{align}
which in this case gives $\E(W_U) = \E(W_U \cond E) (1+o(1))$ and thus
\begin{align} \label{EUE}
 \E(W_U \cond E) &= \E(W_U) (1+o(1)) .
\end{align}

Also by the law of total expectation,
\begin{align} \label{EX}
\E(W_X)
 &= \E(W_X \cond E) \Pr(E) + \E(W_X \cond \Ebar) \Pr(\Ebar) .
\end{align}

As previously established, under event $E$, $W_X = W_U (1+o(1))$.
It follows that
\begin{align} \label{EXE}
\E(W_X \cond E)
 = \E(W_U \cond E) \, (1+o(1))
 = \E(W_U) \, (1+o(1)),
\end{align}
where in the second equality we have used \eqref{EUE}.
In the event $\Ebar$, $W_X$ lies between 0
and the weight of a prescribed copy of $H$
(say, on vertices $1,\ldots,v$, in that order).
We may test for event $E$ by revealing edge weights up to $w_0$ in model $U$,
so for each edge
we know either the exact weight (at most $w_0$), or know that the weight is $\geq w_0$.
In the $X$ model, correspondingly, in the prescribed copy we know the
edge weights up to $x_0 = F^{-1}(w_0)$, or that the weight is $\geq x_0$.
The expected weight of each edge is larger in the case that it is known to be $\geq x_0$
and, even making this pessimistic assumption for every edge in the prescribed copy, we have
\begin{align} \label{EXEbar}
 0 < \E(W_X \cond \Ebar) \leq m \cdot \E(X \cond X>x_0)) = O(1) .
\end{align}
(The conditional expectation $\E(X \cond X>x_0)$ cannot be infinite,
as then by the law of total expectation the expectation of $X$ itself
would be infinite, contradicting our hypothesis.)

Substituting \eqref{EXE} and \eqref{EXEbar} into \eqref{EX} gives
\begin{align}
\E(W_X)
 &= \E(W_U) \, (1+o(1)) + O(1) \exp(-n^{\Omega(1)})
 \sim \E(W_U) .
\end{align}
This completes the proof of the extended expectation result.

\section{Sample Results and Discussion} \label{sample}
In this section we discuss the quality of the results provided by Theorem \ref{th1x},
the lower and upper bounds
--- call them respectively $\Fm(w)$ and $\Fp(w)$ --- on the CDF $F(w)$ of the
weight $W$ of a minimum-weight $k$-clique in
a randomly edge-weighted complete graph of order $n$.

\begin{figure}[htbp]
\begin{center}
\ifpdf{%
\includegraphics[width=3in]{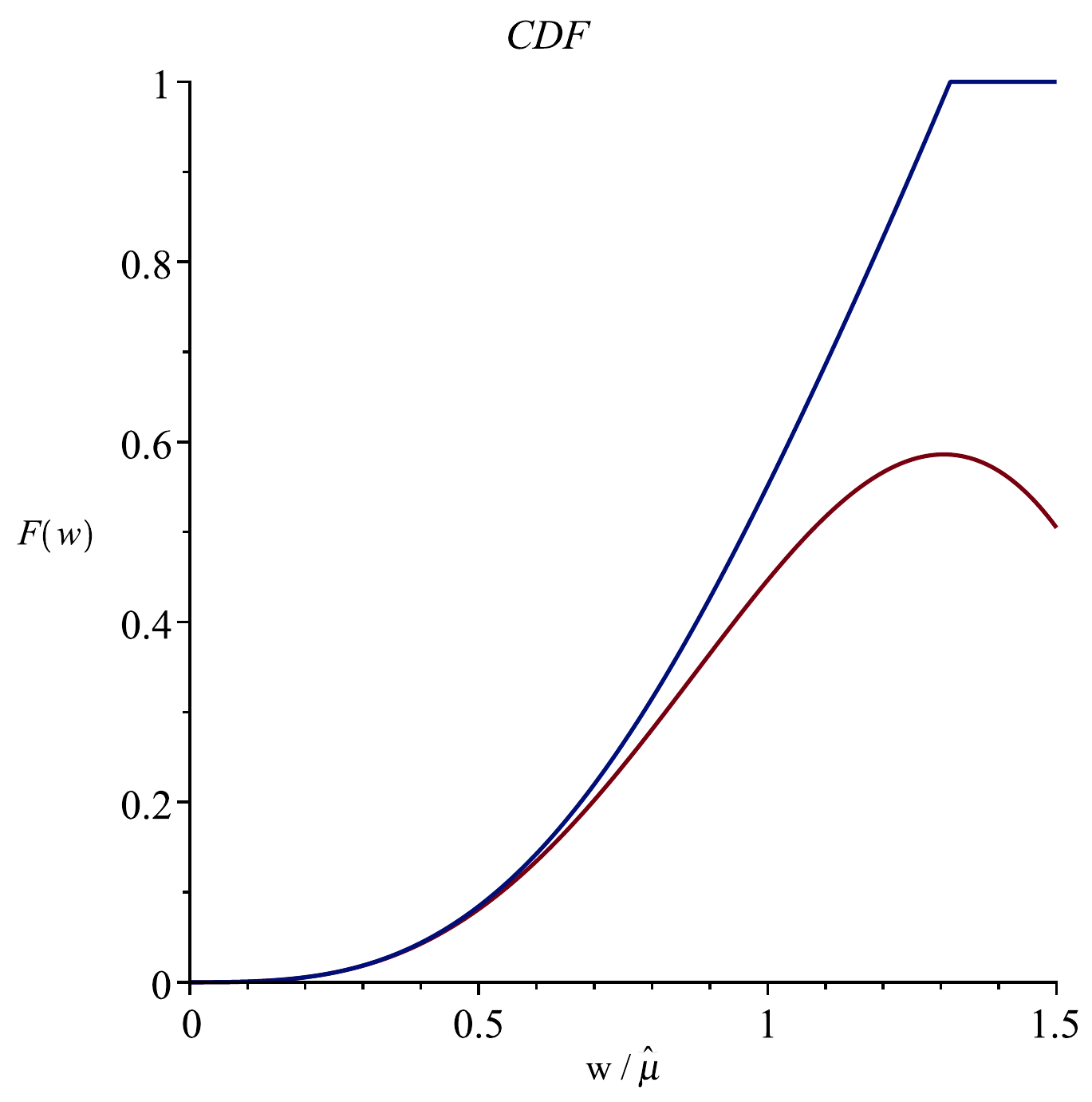}
}\else{%
}\fi
\hspace*{1cm}
\ifpdf{%
\includegraphics[width=3in]{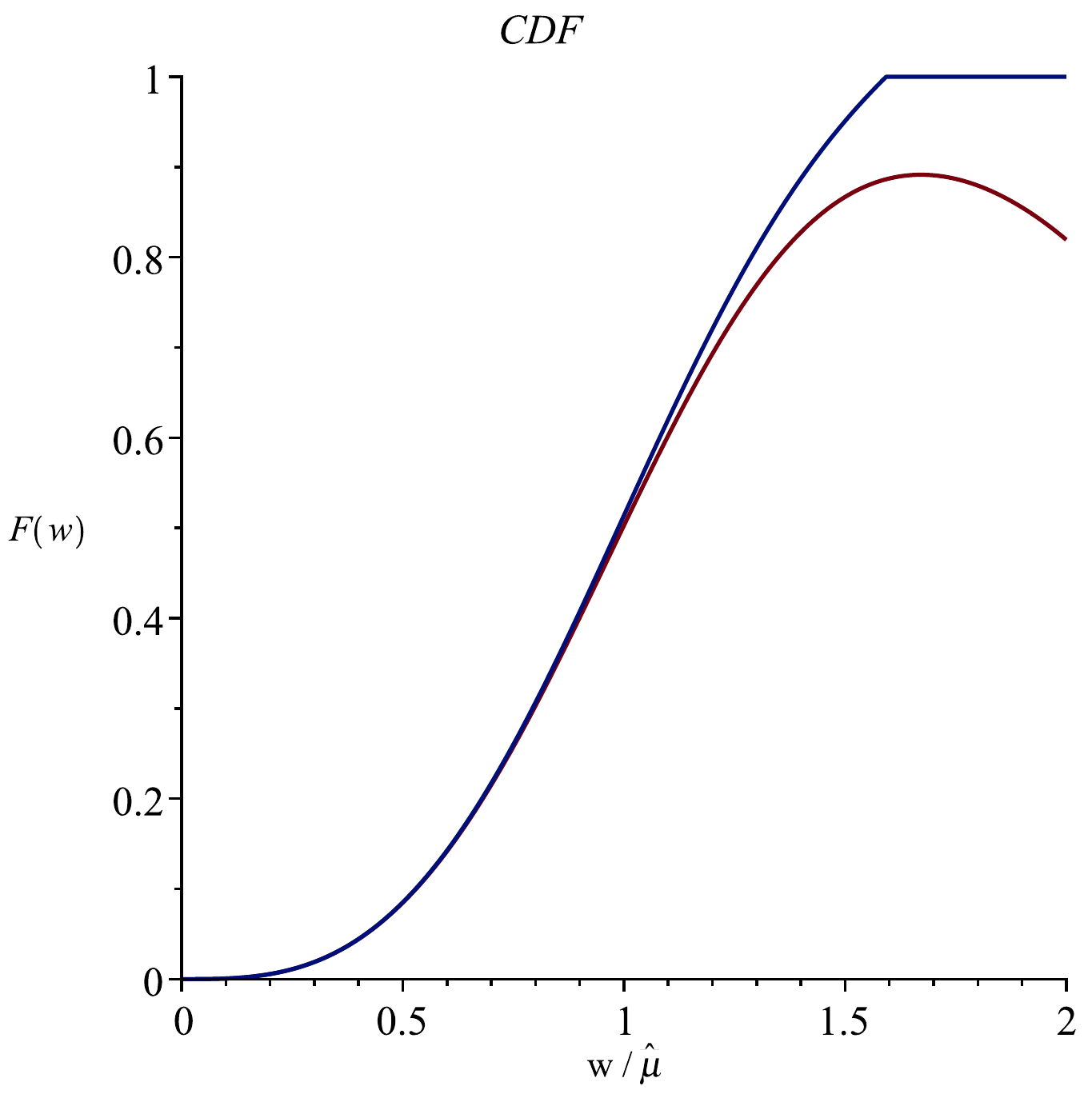}
}\else{%
}\fi
\end{center}
\caption{%
CDFs for
$k=3$, $n=100$ (left) and
$k=3$, $n=1,000$ (right).
}
\label{k3n1000}
\end{figure}

Figure \ref{k3n1000} shows $\Fm$ and $\Fp$
for $k=3$, $n=100$ (left)
and $k=3$, $n=1,000$ (right).
The vertical axis indicates cumulative probability;
the horizontal axis indicates $w$ and is given in units of the estimated mean $\muhat$
given by \eqref{0H} of Theorem \ref{thmH}.
Here $m=\binom k 2$, $d=m/k=(k-1)/2$, and $a=k!$,
and by \eqref{0H} of Theorem \ref{thmH},
$$
\E(W) \sim \muhat = \frac1{n^{2/(k-1)}}
   \frac{\left( { \binom k 2}! k! \right)^{1/\binom k 2}}{\binom k 2}
   \: \Gamma\left( \frac1{\tbinom k 2} \right)
 .
$$
For $k=3$ this gives
$\E(W)
 \sim \muhat = \frac1n \frac{36^{1/3}}{3} \Gamma(1/3)
 = (1.2878\ldots) \, n^{-1}
$.

Looking at Figure \ref{k3n1000}, for $k=3$, $n=100$
we are getting good estimates in the lower tail,
mediocre estimates for values of $w$ near the (estimated) mean,
and poor estimates in the upper tail.
For $k=3$, $n=1,000$ we get good results in the lower tail and through the mean,
but still poor results in the upper tail.

The probability bounds from the theorem can be less than 0 or greater than 1.
In the following discussion and in Figure \ref{k3n1000}, we have truncated both bounds to the range $[0,1]$,
in particular capping $\Fp$ at 1.
Also, because the error terms increase with $W$,
$\Fm$ is not monotone increasing.
In the following discussion we artificially force it to be (weakly) increasing,
by replacing $\Fm(w)$ with $\max\set{w' \leq w \colon \Fm(w') }$.
(To show the nature of the calculated bounds, though,
this was not done in Figure \ref{k3n1000}.)

\begin{table}
\begin{center}
\begin{tabular}{|rr||l|lll|l|l|}
\hline
$k$ & $n$ & 0.05 & $\muhat$ & LB & UB & 0.95 & max gap
\\ \hline
\tabrow {3}{100}{0.04862}{0.02949}{0.44556}{0.55215}{---}{0.41804}
\\ \tabrow {3}{1,000}{0.04986}{0.00295}{0.50278}{0.51387}{---}{0.11282}
\\ \tabrow {3}{10,000}{0.04999}{0.00029}{0.50871}{0.50983}{0.96326}{0.02310}
\\ \tabrow {3}{100,000}{0.05}{0.00003}{0.50931}{0.50942}{0.95124}{0.00403}
\\ \hline
\tabrow {4}{100}{0.04470}{0.21895}{0.31442}{0.58856}{---}{0.66515}
\\ \tabrow {4}{1,000}{0.04947}{0.04717}{0.45684}{0.48193}{---}{0.19120}
\\ \tabrow {4}{10,000}{0.04995}{0.01016}{0.47005}{0.47236}{0.98093}{0.03674}
\\ \tabrow {4}{100,000}{0.05}{0.00219}{0.47127}{0.47149}{0.95236}{0.00594}
\\ \hline
\tabrow {10}{3,000,000}{0.04998}{0.88550}{0.43948}{0.43978}{0.95220}{0.00409}
\\ \tabrow {10}{10,000,000}{0.05}{0.67763}{0.44100}{0.44109}{0.95077}{0.00131}
 \\ \hline
\end{tabular}
\caption{Measurements of the quality of lower and upper bounds on the CDF for various
values of $k$ and $n$.
The column ``0.05'' gives the value of $\Fm(w)$ at the $w$ where $\Fp(w)=0.05$.
The next three columns give the estimated mean $\muhat$ of $W$
along with $\Fm(\muhat)$ and $\Fp(\muhat)$.
The column ``0.95'' gives $\Fp(w)$ where $\Fm(w)=0.95$.
The column ``max gap'' is the largest difference, over all $w$, of $\Fp(w)-\Fm(w)$;
in all cases, it was equal to the gap between the maximum of $\Fm$ and 1.
}
\label{stats}
\end{center}
\end{table}

These preliminary observations on Figure \ref{k3n1000} suggest a few measures of interest,
compiled in Table \ref{stats}.
Let us explain the table and the results observed.

\textbf{Lower tail tests}
The most natural application of our results is to perform lower tail tests.
It is easy to imagine contexts which would result in smaller-weight cliques
than i.i.d.\ edge weights would produce,
for example social networks in which if there is an affinity (modeled as a small weight) between A and B,
and an affinity between B and C, then there is likely also to be an affinity between A and C.

If we wish to show that values of $W$ as small as one observed
occur with probability less than (say) 5\% under the null hypothesis
(that weights are i.i.d.\ uniform $(0,1)$ random variables),
then that observation must be at or below the point $w$ where
$\Fp(w)=0.05$.
If at this point the lower bound is, say, $\Fm(w)=0.02$,
and if the latter happens to be the truth (if $\Fm$ rather than $\Fp$ is a good approximation to $F$ here),
then we require
an observation at the 2\% level
to demonstrate significance at the 5\% level,
and this may prevent our doing so.

We therefore take as a measure of lower-tail performance
the value of $\Fm(w)$ at the point $w$ where $\Fp(w)=0.05$,
that is, $\Fm(w)|_{\Fp(w)=0.05}$.
Show in the table column ``$0.05$'',
if this value is close to $0.05$ our estimates have given up little,
and this is seen largely to be the case throughout the table.

Of course we may be interested in one-tail significance tests at confidence levels $\al < 0.05$
($0.05$ being the largest threshold in common use).
In this case we would hope for a small gap $(\Fp(w)-\Fm(w))|_{\Fp(w)=\al}$
or small ratio $\Fp(w)/\Fm(w)|_{\Fp(w)=\al}$.
Experimentally, both of these measures appear to be increasing functions of $\al$
(i.e., decreasing as $\al$ decreases),
and thus the high quality of our bounds at $\al=0.05$
implies the same for any $\al \leq 0.05$.
Our bounds thus appear to be quite useful for lower tail tests.

\textbf{Mid-range values}
It is natural to check how good our probability bounds are for typical values of $W$.
Taking the estimated mean $\muhat$ of $W$ to stand in for a typical value,
the table reports $\muhat$ and the lower and upper bounds $\Fm(\muhat)$ and $\Fp(\muhat)$
on the CDF $F(w)$ at this point.
It can be seen that the quality of this mid-range estimate is poor for $k=3$, $n=100$,
where the gap $\Fp(\muhat)-\Fm(\muhat)$ is above $0.1$,
but considerably better for $n=1,000$ and $n=10,000$,
where the gap is only about $0.01$ or $0.001$ respectively.

\textbf{Upper tail estimates}
Our results might also be applied to perform upper tail tests.
This would be appropriate for contexts that would produce larger-weight cliques
than i.i.d.\ edge weights would produce,
perhaps a social ``enmity'' network in which if there is enmity (modeled as a small weight) between A and B,
and enmity between B and C, then (on the basis that ``the enemy of my enemy is my friend'')
there is likely to be less enmity (larger weight) between A and C.

Our first measure here is the obvious analogue of the lower-tail one:
the value of $\Fp$ at the point where $\Fm$ is 0.95,
$\Fp(w)|_{\Fm(w)=0.95}$.
This is shown in the table in the column ``0.95''.
In many cases $\Fm$ never even reaches $0.95$, this measure is undefined,
and the implication is that it is impossible to establish upper-tail significance with our method
even at the 5\% confidence level.
However, once $n$ is large enough that the measure is defined,
larger values of $n$ quickly lead to a small gap $(Fp(w)-\Fm(w))|_{\Fp(w)=0.95}$:
if we can in principle establish upper-tail significance,
we can often do so fairly efficiently.

\sloppypar{A second measure relevant here is the maximum gap between our lower and upper bounds,
$\max_{w >0}(\Fp(w)-\Fm(w))$.
Typically the maximum gap is achieved at the smallest point where $\Fp(w)=1$,
and thus if the gap is larger than 0.05
(as for instance for $k=3$, $n=1,000$, for which the maximum gap is around $0.11$)
an upper-tail significance at the 5\% level cannot possibly be established.
We chose this measure rather than the maximum of $\Fm$ because this has the stronger interpretation
that our probability estimates are this accurate across the range:
for any observed $W$, we can report lower and upper bounds on the corresponding CDF value
(under the null hypothesis) no further apart than this gap.
}

\textbf{Parameter values and potential improvements}
For $k=3$, values of $n$ as small as 100 give good estimates in the lower tail,
and modestly good ones for typical values of $W$,
but no upper tail results.
By $n=10,000$, results are good across the range, with a maximum gap of around $0.02$;
with $n=100,000$ this decreases to $0.004$.

For $k=4$ there is a similar pattern, with only slightly less sharp results for the same $n$.

For $k=10$ the picture is significantly different.
Recall that our methods restrict us to estimating $F(w)$ for $w \leq 1$
and here that leaves us hopelessly far into the left tail.
With $k=10$ and $n=100,000$,
the estimated mean of $W$ is $\muhat = 1.8856$, 
while $\Fp(1)$ is less than $10^{-12}$.
For $n=1,000,000$, with $w\leq 1$ we still cannot access the estimated mean $\muhat=  1.13036$,
and $\Fp(1) \approx 0.00231$:
our methods would be useful for observations anomalously small at the 0.1\% confidence level
(to name a standard value near 0.00231),
but nothing much above that.

However, with $k=10$, $W$ is concentrated near $\muhat$,
and thus, \emph{once $n$ is large enough that $\muhat$ falls below 1},
our methods give good results across the range,
as shown in the table.
So, for $k=10$, the problem we observe is with the \emph{range of validity} of our estimates
rather than their quality.
If Theorem \ref{th1x} is extended as outlined in Section \ref{unifx},
the results might well be adequately tight.

One other weakness of our methods is that the lower bound $\Fm$ falls significantly short of 1
for $k=3$ and $k=4$ with $n=100$ and $n=1,000$, revealed in the table's large ``maximum gap'' measures.
It might be possible to improve this by applying the method used in proving
Claim \ref{PrExp},
but we have not attempted this.

\section{Conclusions}
The object studied in this work is the distribution of a minimum-weight clique,
or copy of a strictly balanced graph $H$,
in a complete graph $G$ with i.i.d.\ edge weights.
Theorem \ref{thmH} provides
asymptotic characterizations of the distribution and its mean
for any strictly balanced graph $H$,
while Theorems \ref{th1} and \ref{th1x} provide explicit (non-asymptotic)
descriptions of the distribution for cliques.

This distribution is a natural object of mathematical study,
but also likely to have practical relevance,
particularly for statistical determination that a given network's weights
are \emph{not} i.i.d.
We look forward to seeing such applications of the work.
Some potential applications would involve networks that are not complete graphs,
and extending our results to such cases seems challenging,
whether by extending to other infinite graph classes of graphs
or by including a particular graph as part of the input.

As presented, our explicit methods are for cliques, the uniform distribution
and clique weights at most 1.
However, as discussed in Section \ref{th1xx}, it is easy to write down calculations
for the uniform distribution and all clique weights, the exponential distribution,
and probably other common distributions.
In doing so we encountered computational challenges,
but these seem surmountable if the incentive is more than just fleshing out a table.
As noted in Section \ref{exactH}, extending explicit results to subgraphs $H$ other than cliques
is straightforward.

The quality of the results from our methods, as discussed in Section \ref{sample},
is largely good, especially when we are interested in lower-tail results
and relatively small cliques, or larger cliques in very large graphs.
To improve the probability bounds in the middle range,
near the mean, would seem to require an approach other than Stein-Chen,
but we have no concrete alternative suggestions.
Improving the bounds in the upper tail might be done, as suggested earlier,
by applying the ideas of Claim \ref{PrExp}.

\end{document}